\documentclass{article}

\usepackage[T1]{fontenc}
\usepackage[utf8]{inputenc}
\usepackage[british]{babel}

\usepackage{amsmath}
\usepackage{amsthm}
\usepackage{amsfonts}
\usepackage{hyperref}
\hypersetup{
	colorlinks = true,
	linkcolor = purple, 
	filecolor = purple, 
	citecolor = purple,
	urlcolor = purple
}
\usepackage{url} 
\usepackage{amssymb}
\usepackage{stackrel}
\usepackage{tikz-cd}
\allowdisplaybreaks

\newtheorem{lemma}{Lemma}
\newtheorem{theorem}[lemma]{Theorem} 
\newtheorem{prop}[lemma]{Proposition}

\theoremstyle{definition}
\newtheorem{rem}[lemma]{Remark}
\newtheorem{nota}[lemma]{Notation}
\newtheorem{defn}[lemma]{Definition}

\newcommand{\Z}{\mathbb{Z}}
\newcommand{\Q}{\mathbb{Q}}
\newcommand{\C}{\mathbb{C}}
\newcommand{\R}{\mathbb{R}}
\newcommand{\A}{\mathbb{A}}
\newcommand{\SH}{\mathcal{SH}}

\DeclareMathOperator{\Sym}{Sym}
\DeclareMathOperator{\Tr}{Tr}
\DeclareMathOperator{\GW}{GW}
\DeclareMathOperator{\sgn}{sgn}
\DeclareMathOperator{\Spec}{Spec}
\DeclareMathOperator{\End}{End}
\DeclareMathOperator{\spn}{sp}

\newcommand{\Var}{\text{Var}}

\title{Quadratic Euler Characteristic of Symmetric Powers of Curves}
\author{Lukas F. Bröring and Anna M. Viergever}
\date{}

\begin{document}
\maketitle

\begin{abstract}
  We compute the quadratic Euler characteristic of the symmetric powers of a
  smooth, projective curve over any field $k$ that is not of characteristic
  two, using the Motivic Gauss-Bonnet Theorem of Levine-Raksit. As an
  application, we show over a
  field of characteristic zero that the power structure on the Grothendieck-Witt ring
  introduced by Pajwani-P\'al computes the compactly supported $\mathbb{A}^1$-Euler
  characteristic of symmetric powers for all curves.
\end{abstract}

\tableofcontents 
\phantomsection
\addcontentsline{toc}{section}{Introduction}

\section*{Introduction}
One can assign a \textit{quadratic Euler characteristic} to a smooth,
projective scheme over a field $k$ that is not of characteristic $2$ using motivic
homotopy theory, introduced by Morel-Voevodsky. Quadratic Euler
characteristics were first introduced by Hoyois \cite{HoyoisQRGLVTF} and are
elements of the Grothendieck-Witt ring $\GW(k)$ of quadratic forms over
$k$. These quadratic forms carry a lot of information: If
$k\subset\mathbb{R}$, Levine \cite{LevineA} proved that the signature of the
quadratic Euler characteristic $\chi(X/k)$ of a smooth projective scheme $X$
over $k$ is equal to the topological Euler characteristic of $X(\mathbb{R})$
and the rank of $\chi(X/k)$ is the topological Euler characteristic of
$X(\mathbb{C})$. In practice however, $\chi(X/k)$ is generally hard to
compute. Quadratic Euler characteristics are often used in the program of
refined enumerative geometry, which aims to obtain ``quadratically enriched''
versions of results in classical enumerative geometry.

Our main theorem is the following.
\begin{theorem}\label{main theorem}
	Let $C$ be a smooth, projective curve of genus $g$ over a field $k$
        that is not of characteristic $2$. Let $n\in\mathbb{Z}_{\geq 1}$.
        Then if $n = 2m$ is even, we have
	\begin{align*}
          \chi(\Sym^{n}C/k)
          &= \sum_{i=0}^m\binom{g}{i}\langle -1\rangle^i +
            \frac{1}{2}\left(\binom{2g-2}{n} -
            \sum_{i=0}^m\binom{g}{i}\right)H\in \GW(k)\\
          &= \frac 12\left(\binom{2g-2}{n} + (-1)^m\binom{g-1}{m}\right)\cdot
            \langle 1\rangle \\&\quad+ \frac 12
            \left(\binom{2g-2}{n}-(-1)^m\binom{g-1}{m}\right)\cdot\langle -1\rangle
	\end{align*}
	and if $n$ is odd, we have
	\[
		\chi(\Sym^{n}C/k) = -\frac 12 \binom{2g-2}{n} H\in \GW(k),
	\]
	where $\langle a\rangle\in\GW(k)$ for $a \in k^\times$ denotes the quadratic
        form $x \mapsto ax^2$ and $H = \langle 1\rangle + \langle -1\rangle\in\GW(k)$ is the hyperbolic
        form.
\end{theorem}

\begin{rem}
  For a smooth, projective curve $C$, the schemes $\Sym^nC$ are again smooth and
  projective. Thus, the quadratic Euler characteristic of $\Sym^nC$ is
  well defined.
\end{rem}

\begin{rem}\label{rem:rank}
  MacDonald \cite[(4.4)]{MacDonaldS} computed the topological Euler
  characteristic of $\Sym^nC$ for a smooth, projective complex curve $C$ to
  be
  \[
    \chi^{top}(\Sym^nC) = (-1)^n\binom{2g-2}{n}.
  \]
  This implies that if $C$ is a
  smooth, projective curve over a
  field $k$ of characteristic zero, we have
  \[
    \operatorname{rank}(\chi(\Sym^nC/k)) = (-1)^n\binom{2g-2}{n},
  \]
  which matches with the rank of the forms in Theorem \ref{main theorem}. Indeed, there exists a smooth, projective curve $C_0$ over a
  subfield $k_0 \subset k$ with the following two properties: First, $k_0$ is
  finitely generated over $\Q$ and, second, the base change $(C_0)_k$ to $k$
  is isomorphic to $C$. Since the rank of the quadratic Euler characteristic
  is invariant under base-change and $k_0$ admits an embedding into $\C$, we get after choosing an embedding $k_0 \to \C$ that
  \begin{align*}
    \operatorname{rank}(\chi(\Sym^nC/k)) &= \operatorname{rank}(\chi(\Sym^nC_0/k_0))\\
    &=
    \operatorname{rank}(\chi(\Sym^n(C_0)_\C/\C))\\
    &= (-1)^n\binom{2g-2}{n}
  \end{align*}
  as desired.
\end{rem}

\begin{rem}\label{remark: vanishing and rank}
  The formula in Theorem \ref{main theorem} evaluates to zero for a smooth,
  projective curve $C$ of genus $g$ whenever $n > 2g-2$. If $C$ has a rational
  point, this can also be seen with the following brief argument: Consider the map
  $\Sym^nC\to \text{Pic}^n(C)$ sending a point on $\Sym^nC$ with
  associated divisor $D$ to $\mathcal{O}_C(D)$. By for example the proof of \cite[Theorem 7.33]{MustataZF}, this map realizes $\Sym^nC$ as
  a projective space bundle over $\text{Pic}^n(C)$ and we can find an
  isomorphism $\text{Pic}^n(C)\cong \text{Pic}^0(C)$. Using the Motivic
  Gauss-Bonnet Theorem, the quadratic Euler characteristic of any
  Abelian variety is zero because its tangent bundle is trivial. Combining
  this with \cite[Proposition 1.4]{LevineA} yields the desired vanishing.
\end{rem}

% Our application
As an application, we extend a result on compatibility of power structures of
Pajwani and P\'al \cite{Pajwani-PalPS}:
Over a field $k$ of characteristic
zero, Arcila-Maya, Bethea, Opie, Wickelgren and Zakharevich
\cite{Wickelgren} extended the quadratic Euler characteristic to a motivic measure
\[\chi_c: K_0(\text{Var}_k)\to \GW(k)\]
called the compactly supported $\mathbb{A}^1$-Euler characteristic. On $K_0(\text{Var}_k)$, Gusein-Zade,
Luengo and Melle-Hernández \cite{Gusein-Zade-Luengo-Melle-HernandezPS} showed that
symmetric powers of schemes define a power structure. In \cite{Pajwani-PalPS},
Pajwani and P\'al constructed a power
structure $a_\ast$ on $\GW(k)$ with the following property: If there exists a power
structure on $\GW(k)$ such that $\chi_c$ respects the power structures, then the power structure on $\GW(k)$ has to agree with
$a_\ast$. It is still an open question,  whether there indeed is a power
structure on $\GW(k)$ such that $\chi_c$ respects the power structures. In
Theorem \ref{theorem power structures respected}, we prove that if $\GW(k)$ is given the power structure $a_*$, then $\chi_c$ respects the power
structures for all curves over a field of characteristic zero.

% Why is our result interesting
At the moment, there is no general formula for the quadratic Euler
characteristic of a quotient scheme. Theorem \ref{main theorem} is an example
result for this
and might be helpful for understanding the general case. Further work
on quadratic Euler characteristics of quotient schemes will appear in the
upcoming PhD thesis of the first named author.

% Strategy of the proof
The proof of Theorem \ref{main theorem} uses an explicit computation in
conjunction with the Motivic
Gauss-Bonnet Theorem proven by Levine and Raksit \cite{LevineGB}, which computes
the quadratic Euler characteristic of a smooth projective
scheme as a composition of a cup product and trace on Hodge cohomology.

% Related works
In related work,
the Motivic Gauss-Bonnet Theorem was also used by Levine, Pepin Lehalleur, and Srinivas \cite{LevineLehalleurSrinivas}
to compute the quadratic Euler characteristic of a smooth hypersurface in
$\mathbb{P}^n$ and in the second named author's computation of the quadratic Euler
characteristic of a smooth same-degree
complete intersection in $\mathbb{P}^n$, see \cite{Viergever}.

Adjacent to our application, Pajwani and P\'al \cite{Pajwani-PalPS} showed that their power structure is
respected for zero dimensional schemes. Work of Pajwani, Rohrbach, and the
second named author
\cite{PajwaniRohrbachViergever} extended this to the class of ``étale linear
varieties'', which is the subring of $K_0(\text{Var}_k)$ generated by the class
of $[\mathbb{A}^1]$ and the classes of the form $[\text{Spec}(L)]$ where $L$ is a finite \'etale algebra over
$k$. Since curves are almost always outside
of the class of ``étale linear varieties'' by \cite[Corollary
6.6]{PajwaniRohrbachViergever}, Theorem \ref{theorem power structures
  respected} yields a proper extension of the class of
varieties for which the power structure is respected.

We would like to highlight that Theorem \ref{main theorem} has been proven independently by Simon Pepin Lehalleur and Lenny Taelman, as part of work in progress on the higher dimensional case.

\subsection*{Acknowledgments}
We would like to extend our sincere thanks to Marc Levine for suggesting this project to us, for helpful
discussions, and for useful comments on an earlier draft of this paper. We are very grateful to Lenny Taelman and Simon Pepin
Lehalleur for their helpful comments on an earlier draft of this paper. We also thank Jens Hornbostel for helpful discussions. We further thank Jesse Pajwani for pointing us to the paper
\cite{McGarraghySP} and for helpful comments on an earlier draft of this
paper. We thank Herman Rohrbach for pointing us to the paper
\cite{HoyoisHKOGA} and for helpful discussions. In addition, we
would like to acknowledge Raymond Cheng and Hind Souly for helpful discussions.

The first named author was supported by the RTG 2553. The second named author was partially supported by the European Research Council (ERC) under the European Union’s Horizon 2020 research and innovation programme under grant agreement No 948066 (ERC-StG RationAlgic).

\subsection*{Notation and Conventions}
Throughout, we let $k$ be a field that is not of characteristic $2$. By a
curve, we mean a one-dimensional, geometrically connected, separated, finite-type scheme over $k$. For a
smooth projective scheme $X$ over $k$, we write $\Omega_X$ for the sheaf of
K\"ahler differentials on $X$. For $i\in\mathbb{Z}_{\geq 0}$, we write
$\Omega_X^i = \Lambda^i\Omega_X$.

For $p,q, p', q' \in \Z_{\ge 0}$, we denote the cup product of $a \in
H^q(X,\Omega^p_X)$  and $b \in H^{q'}(X,\Omega^{p'}_X)$ by $ab := a\cdot b :=
a \cup b \in H^{q+q'}(X,\Omega^{p+p'}_X)$.
This cup-product is graded commutative if we
consider $a \in H^q(X,\Omega^p_X)$ to be of degree $q-p$. That is,  we have for
$a \in H^q(X,\Omega^p_X)$ and $b \in H^{q'}(X,\Omega^{p'}_X)$ that $ab =
(-1)^{(q-p)(q'-p')} ba$.

We let $K_0(\text{Var}_k)$ be the Grothendieck group of varieties over
$k$. This is the free Abelian group generated by isomorphism classes $[X]$ of
varieties $X/k$ modulo the relation that for $Z\subset X$ a closed
subvariety we have $[X] = [Z] + [X\setminus Z]$. The product rule $[X][Y] =
[X\times_kY]$ induces a ring structure on~$K_0(\text{Var}_k)$.

For a quasi-compact and quasi-separated scheme $X$, we denote its stable
motivic homotopy category as described by Hoyois \cite{HoyoisSOEMHT} by $\SH(X)$. We denote the
unit object with respect to the smash product in $\SH(X)$ by $1_X\in \SH(X)$.

\section{Quadratic Euler Characteristic}
The quadratic Euler characteristic, which was first studied by Hoyois
\cite{HoyoisQRGLVTF}, is a refinement of the
topological Euler characteristic. We briefly recall its construction. 

To $k$, we can associate the stable motivic
homotopy category $\SH(k)$ constructed by Morel-Voevodsky \cite{MorelATF}, which is a symmetric monoidal category. See
also Hoyois \cite{HoyoisSOEMHT} for an introduction to $\SH(k)$. Given a
smooth projective scheme $X/k$ with structure map $p_X: X\to \Spec(k)$,
the object $(p_X)_*1_X\in \SH(k)$ is strongly dualizable by for example \cite[Theorem 5.22]{HoyoisSOEMHT}. One can now apply the categorical Euler characteristic of Dold-Puppe \cite{dold80-categoricalEC} to $(p_X)_*1_X$, which yields an endomorphism of the unit of $\SH(k)$. Morel has proven that $\End_{\SH(k)}(1_k)\cong \GW(k)$, and thus we obtain the quadratic Euler characteristic $\chi(X/k)\in \GW(k)$ as an element of the Grothendieck-Witt ring of quadratic forms over $k$,
which we introduce below.

\begin{defn}
  The \emph{Grothendieck-Witt ring $\GW(k)$} over $k$ is
  defined as the group completion of the monoid of isometry classes of
  non-degenerate quadratic forms over $k$ with respect to taking orthogonal
  direct sums.

  The group $\GW(k)$ inherits a ring structure from the tensor
  product of quadratic forms.
\end{defn}  

\begin{rem}
  Over a field of characteristic not $2$, a non-degenerate quadratic form is
  the same as a non-degenerate symmetric bilinear form, and elementary linear
  algebra computations yield the following presentation. The
  quadratic forms $\langle a \rangle:x \mapsto ax^2$ for $x \in k^\times$
  generate $\GW(k)$ and they are subject to the following relations for
  $a, b\in k^\times$:
  \begin{itemize}
  \item $\langle a\rangle \cdot \langle b\rangle = \langle a b \rangle$,
  \item
    $\langle a\rangle + \langle b\rangle = \langle a+b\rangle + \langle
    ab(a+b)\rangle$, whenever $a+ b \in k^\times$,
  \item $\langle ab^2\rangle = \langle a\rangle$, and
  \item
    $\langle a\rangle + \langle -a\rangle = \langle 1\rangle + \langle
    -1\rangle =: H$.
  \end{itemize}
  The quadratic form $H = \langle 1\rangle + \langle -1\rangle$ is called the
  \emph{hyperbolic form}.

  Witt \cite[Section 1]{WittQF} first described this presentation. In the form
  presented here, it is
  \cite[Lemma 3.9]{MorelATF}. There, the statement is deduced from the
  statement for Witt rings, proven by Milnor and Husemoller
  \cite[Lemma (1.1) in Chapter 4]{MilnorSBF}.
\end{rem}

See \cite[Section
2]{LevineA} for a detailed exposition of the quadratic
Euler characteristic and its basic
properties.

\begin{rem}
  In \cite[Remark 2.3]{LevineA}, Levine shows that for $k \subset \R$ and $X$
  a smooth, projective scheme over $k$, the rank of $\chi(X/k)$ agrees with
  the topological Euler characteristic of $X(\mathbb{C})$ and the signature
  agrees with the topological Euler characteristic of $X(\R)$. By a theorem of Saito \cite[Theorem 2]{Saito}, the discriminant of $\chi(X/k)$ can also be interpreted in terms of the determinant of $\ell$-adic cohomology for any $\ell$ which is coprime to the characteristic of $k$. See  \cite[Theorem 2.22]{Pajwani-PalYZ} and the surrounding text, for a more detailed discussion. 
\end{rem}

The Motivic Gauss-Bonnet Theorem of Levine-Raksit \cite{LevineGB} provides a
way to compute the quadratic Euler characteristic using Hodge
cohomology. A more general version has been proved by
Déglise-Jin-Khan \cite{DegliseFCMHT}. We do not state the theorem in full
generality here, but rather a simplified version of one of its corollaries,
which is enough to prove Theorem \ref{main theorem}.

\begin{theorem}[Motivic Gauss-Bonnet, {\cite[Corollary 8.7]{LevineGB}}]
  Let $X$ be a smooth, projective scheme over $k$.
  \begin{itemize}
  \item If $\dim X$ is odd, then $\chi(X/k) = m\cdot H$ for some $m \in \Z$.
  \item If $\dim X = 2n$ is even, then $\chi (X/k) = Q + m\cdot H$ for some $m \in
    Z$ where $Q$ is the quadratic form given by the composition
    \[
      H^n(X,\Omega^n_X)\times H^n(X,\Omega^n_X) \xrightarrow{\cup}
      H^{2n}(X,\Omega_X^{2n}) \xrightarrow{\Tr} k.
    \]
    Here, $\cup$ denotes the cup product and $\Tr$ the trace map from Serre
    duality.
  \end{itemize}
\end{theorem}

Arcila-Maya, Bethea, Opie, Wickelgren, and Zakharevich \cite{Wickelgren} constructed a motivic
measure $\chi_c\colon K_0(\Var_k)\to \GW(k)$, called the compactly
supported $\A^1$-Euler characteristic, which extends the quadratic Euler
characteristic of smooth, projective schemes
to varieties over $k$ when $k$ has characteristic 0.  Levine, Pepin Lehalleur and
Srinivas \cite[Section
5.1]{LevineLehalleurSrinivas} give a definition of $\chi_c$ which works for any field which is not of characteristic two. We briefly recall the construction here. 

If $X$ is a variety over a field $k$ of characteristic zero with structure map
$p_X: X\to\text{Spec}(k)$, one can show using resolution of singularities that
$(p_X)_!(1_X)\in \SH(k)$ is a strongly dualizable object. The categorical
Euler characteristic of $(p_X)_!(1_X)$ then computes the
compactly supported
$\mathbb{A}^1$-Euler characteristic $\chi_c(X/k)\in \GW(k)$. 

If $X$ is a variety over a field $k$ of
characteristic $p>0$ with structure map $p_X: X\to\text{Spec}(k)$, one can
show using results of Riou that $(p_X)_!(1_X)\in \SH(k)$ is a strongly
dualizable object after passing to the perfect completion $k^{perf}$ of $k$ and inverting
$p$. We have that $\GW(k)\cong\GW(k^{perf})$ where $k$ is the perfect completion of $k$, see \cite[Remark 2.1.2]{LevineA}. The categorical Euler characteristic of $(p_X)_!(1_X)$ then yields the
Euler characteristic $\chi_c(X/k)\in \GW(k)[\frac{1}{p}]$. One then shows that
$\chi_c(X/k)$ lies in the image of the injective map
$\GW(k)\to\GW(k)[\frac{1}{p}]$, so that one can define the compactly supported
$\mathbb{A}^1$-Euler characteristic as an element $\chi_c(X/k)\in \GW(k)$. 

\section{Proof of Theorem \ref{main theorem}}
Let $C$ be a smooth, projective curve of genus $g$ over $k$ and let
$n\in\mathbb{Z}_{\ge 1}$. In order to prove Theorem \ref{main theorem}, we first reduce to the case that $k$ has
characteristic zero.

For a discrete valuation ring $A$ with residue field $k$, field of
fractions $K$, and chosen uniformiser $t$, Morel \cite[Lemma 3.16, discussion on p. 57]{MorelATF}
describes a specialisation homomorphism $\spn_t\colon \GW(K) \to \GW(k)$
sending the element $\langle ut^\nu\rangle \in \GW(K)$ with $u \in A^\times$ and $\nu
\in \Z$ to the element $\langle \bar u\rangle \in \GW(k)$. Here $\bar u \in k$
denotes the residue class of $u$. 

The following result is probably known in the study of quadratic conductor formulas, but we were unable to find a proof in the literature. We include a proof here for the reader's convenience.

\begin{lemma}
  \label{lemma:euler-characteristic-deformation}
  Let $A$ be a complete discrete valuation ring with perfect residue field $k$
  and perfect field of fractions $K$. Choose a uniformiser $t \in A$ of
  $A$. Let $X$ be a smooth, projective scheme over $\Spec A$ and
  denote by $X_k$ and $X_K$ its $k$- and $K$-fibre respectively. Then we have
  \[
    \chi(X_k/k) = \spn_t(\chi(X_K/K)).
  \]
  Note that this is independent of the choice of uniformiser $t$.

  \begin{proof}
    We use the six functor formalism on $\SH(-)$ to prove the
    lemma, see for example \cite{HoyoisSOEMHT}. Let $p_A: X\to\Spec(A)$,  $p_k: X_k\to \Spec(k)$, and $p_K: X_K\to \Spec(K)$ be the structure maps of $X$ over $A$, $X_k$ over $k$, and $X_K$ over $K$, respectively.
    Since $X$ is smooth and projective, $(p_A)_!(1_{X}) =
    (p_A)_\ast(1_{X})$ is strongly dualisable by \cite[Theorem
    5.22]{HoyoisSOEMHT}. Thus, the categorical Euler
    characteristic of $(p_A)_!(1_{X})$ in $\SH(A)$ is defined and we denote
    it by
    \[
      \chi^{\text{cat}}(X) := \chi^{\text{cat}}((p_A)_!(1_{X})) \in \End_{\SH(A)}(1_A)
    \]
    and similarly $\chi^{\text{cat}}(X_k) :=
    \chi^{\text{cat}}((p_k)_!(1_{X_k}))\in \End_{\SH(k)}(1_k)$ and $\chi^{\text{cat}}(X_K) :=
    \chi^{\text{cat}}((p_K)_!(1_{X_K}))\in \End_{\SH(K)}(1_K)$.

    Consider the following commutative diagram
    \[
      \begin{tikzcd}
        \SH(X_k)\ar[d, "(p_k)_!"] & \SH(X)\ar[d, "(p_A)_!"]\ar[l,
        "i^\ast"]\ar[r, "j^\ast"] & \SH(X_K)\ar[d, "(p_K)_!"]\\
        \SH(k) & \SH(A)\ar[l, "i^\ast"]\ar[r, "j^\ast"] & \SH(K).
      \end{tikzcd}
    \]
    Since $i^\ast$ and $j^\ast$ are symmetric monoidal, we get using
    \cite[Remark 2.3.1]{LevineA}
    \begin{align*}
      \chi^{\text{cat}}(X_k) &= \chi^{\text{cat}}((p_k)_!1_{X_k})\\ &=
      \chi^{\text{cat}}((p_k)_!i^\ast1_{X}) \\&=
      \chi^{\text{cat}}(i^\ast(p_A)_!1_{X}) \\&=
      i^\ast(\chi^{\text{cat}}((p_A)_!1_{X})) \\&= i^\ast(\chi^{\text{cat}}(X))
      \in \End_{\SH(k)}(1_k)
    \end{align*}
    and similarly $\chi^{\text{cat}}(X_K) = j^\ast(\chi^{\text{cat}}(X))\in \End_{\SH(K)}(1_K)$.

    Hornbostel \cite[Theorem 5.5]{HornbostelARHKW} constructed a
    $\mathbb{P}^1$-ring spectrum $\mathrm{KQ}_S \in \SH(S)$ for $S$ a regular
    scheme with $2 \in \mathcal{O}_S(S)^\times$ representing Hermitian K-theory.
    Using the unit map $u\colon 1_{-} \to
    \mathrm{KQ}_{-}$ of Hermitian K-theory, we can
    map each of the above categorical Euler characteristics to the $(0,0)$-th stable homotopy group of
    $\mathrm{KQ}_{-}$, which is the Grothendieck-Witt ring by for example \cite[Introduction]{HornbostelARHKW}. That is, we
    get $u\chi^{\text{cat}}(X_k) = i^\ast (u\chi^{\text{cat}}(X))\in \mathrm{KQ}^{0,0}_k =
    \GW(k)$ and $u\chi^{\text{cat}}(X_K) = j^\ast (u\chi^{\text{cat}}(X))\in KQ^{0,0}_K =
    \GW(K)$, where we use Hoyois, Jelisiejew, Nardin, and Yakerson's result
    \cite[Lemma 7.5]{HoyoisHKOGA} that Hermitian K-theory is stable under base
    change. The map $u$ is an
    isomorphism in $\SH(k)$, $\SH(A)$ and $\SH(K)$ by an argument of
    Bachmann-Hoyois \cite[Theorem 10.12]{BachmannHoyoisNMHT} and in $\SH(k)$
    and $\SH(K)$ the isomorphisms induced by the unit map agree with Morel's
    identification of the endomorphisms of
    the unit in $\SH(-)$ with the Grothendieck-Witt ring \cite[Theorem 6.4.1
    and Remark 6.4.2]{MorelIAHT}. Thus, we have $\chi(X_k/k) =
    \chi_c(X_k/k) = u\chi^{\text{cat}}(X_k)$ and $\chi(X_K/K) = \chi_c(X_K/K) = u\chi^{\text{cat}}(X_K)$.

    We are now in the following situation:
    \[
      \begin{tikzcd}[column sep=small]
        \GW(k) \ar[d, phantom, "{\rotatebox[origin=c]{90}{$\in$}}"] & \ar[l, "i^\ast"'] 
        \GW(A) \ar[r, "j^\ast"]\ar[d, phantom, "{\rotatebox[origin=c]{90}{$\in$}}"] &
        \GW(K)\ar[d, phantom, "{\rotatebox[origin=c]{90}{$\in$}}"]\\
        \chi(X_k/k) & \ar[l, mapsto] u\chi^{\text{cat}}(X) \ar[r, mapsto] &
        \chi(X_K/K).
      \end{tikzcd}
    \]
    The quadratic form $u\chi^{\text{cat}}(X)$ can be expressed as $\sum_i
    \langle u_i\rangle - \sum_j\langle u'_j\rangle$ for some $u_i, u_j' \in A^\times$ by \cite[Chapter 1, Theorem
    6.4]{ScharlauQHF}. Thus by applying $j^\ast$, we get $\chi(X_K/K) =
    \sum_i\langle
    u_i\rangle - \sum_j\langle u'_j\rangle$, and by applying $i^\ast$, we get $\chi(X_k/k) = \sum_i\langle
    \bar u_i\rangle - \sum_j\langle \bar u'_j\rangle$. In particular, we get by the definition of the
    specialisation morphism that $\chi(X_k/k) = \spn_t(\chi(X_K/K))$.
  \end{proof}
\end{lemma}

\begin{prop}
  Assume Theorem \ref{main theorem} holds over all fields of characteristic
  zero. Then Theorem \ref{main theorem} holds over all fields of characteristic
  not two.

  \begin{proof}
    Assume that $k$ has positive characteristic. Using \cite[Remark
    2.1.2]{LevineA}, we can assume without loss of generality that $k$ is a
    perfect field. By for example
    \cite[Theorem 3 in Chapter II]{SerreLF}, we can find a complete
    discrete valuation ring $A$ with residue field $k$ and field of fractions
    $K$ of characteristic zero. Now \cite[Exposé III, Corollaire 7.4]{GrothendieckSGA1} proves that
    we can lift $C$ to a smooth, projective curve $\tilde C$ over $A$; that is
    we get the following diagram of fibre squares:
    \[
      \begin{tikzcd}
        C\ar[d, "p_k"] \ar[r, "i"]\ar[dr, phantom, "\ulcorner", very near start] & \tilde C \ar[d, "p_A"] & \tilde C_K\ar[l,
        "j"]\ar[d, "p_K"]\ar[dl, phantom, "\urcorner", very near start]\\
        \Spec k \ar[r, "i"] & \Spec A & \Spec K\ar[l, "j"].
      \end{tikzcd}
    \]
    Since taking $\Sym^nC$ is compatible with taking the special and the
    generic fibre, we get the same diagram of fibre squares with $\Sym^nC$,
    $\Sym^n\tilde C$, and $\Sym^n\tilde C_K$ instead of $C$, $\tilde C$, and
    $\tilde C_K$, respectively. Furthermore, $\Sym^nC, \Sym^n\tilde C$ and
    $\Sym^n\tilde C_K$ are smooth and projective, since $\tilde C$ is a smooth,
    projective curve.
    
    By assumption, we know that $\chi(\Sym^n\tilde C_K/K)$ is given by the
    formula in Theorem~\ref{main theorem}. For any choice of uniformiser $t
    \in A$, the specialisation map $\spn_t\colon \GW(K) \to \GW(k)$ maps this
    formula to the same formula in $\GW(k)$. Thus, Lemma
    \ref{lemma:euler-characteristic-deformation} yields the claim.
  \end{proof}
\end{prop}

\begin{proof}[Proof of Theorem \ref{main theorem} in case $n$ is odd and $k$
  has characteristic zero] 
  Since $\Sym^{n}C$ has odd dimension $n$,
  its quadratic Euler characteristic is a multiple of hyperbolic
  forms by the Motivic Gauss-Bonnet Theorem, and thus it is completely
  determined by its rank. Therefore, Remark \ref{rem:rank} implies
  \[
    \chi(\Sym^nC/k) = -\frac 12\binom{2g-2}{n}\cdot H
  \]
  as desired.
\end{proof}

From now on, assume that $n=2m$ is even and that $k$ has characteristic
zero. The only reason for the assumption that $k$ has characteristic zero
is to ensure that $n!$ is invertible in $k$. We will compute $\chi(\Sym^{n}C/k)$ using the Motivic Gauss-Bonnet Theorem. MacDonald \cite{MacDonaldS} has computed the topological Euler characteristic of the symmetric product of a curve over the complex numbers using an explicit basis computation. We follow the same strategy. 

\begin{nota}
  \label{nota:basis}
   Let $\beta^\vee\in H^0(C,\mathcal{O}_C)$ be the unit element with respect
   to the cup product and let $\beta\in H^1(C,\Omega_C^1)$ be the dual element
   with respect to Serre duality. In particular, the trace map sends $\beta$
   to $1$. Let $\alpha_1,\dots, \alpha_g$ be a basis for $H^0(C,\Omega_C)$ and
   let
   $\alpha_1^\vee, \dots, \alpha_g^\vee$ be the dual basis for
   $H^1(C,\mathcal{O}_C)$ with respect to Serre duality.

  Let $\pi_j\colon C^n \to C$ be the projection on the $j$-th
  component. We write $\alpha_i^{(j)} := \pi_j^\ast \alpha_i$ and
  similarly $\alpha_i^{\vee, (j)} := \pi_j^\ast\alpha_i^{\vee}$,
  $\beta^{(j)} := \pi_j^\ast\beta$, and $\beta^{\vee, (j)} := \pi_j^\ast\beta^\vee$.
\end{nota}

\begin{rem}
  \label{rem-multiplication-on-curve}
  We have the following multiplication rules for the above generators. Let
  $i,j\in\{1,\dots, g\}$ be such that $i\neq j$. Then
  \[
    \begin{matrix} 
      \beta\cdot\beta = 0 & & \beta^\vee \cdot \alpha_i = \alpha_i & & \alpha_i\cdot \alpha_i = 0  & & \alpha_i^\vee\cdot\alpha_i^\vee
      = 0\\
      \beta\cdot\alpha_i = 0 & & \beta^\vee\cdot\beta^\vee = \beta^\vee & & \alpha_i\cdot\alpha_j= 0 & & \alpha_i^\vee\cdot \alpha_j^\vee = 0\\
      \beta\cdot \beta^\vee = \beta& & \beta^\vee \cdot \alpha^\vee_i = \alpha^\vee_i&                                                               &\alpha_i^\vee\cdot\alpha_i = \beta & & \\
      \beta\cdot \alpha_i^\vee = 0 & &  & &\alpha_i^\vee\cdot\alpha_j =0 & & 
    \end{matrix}
  \]
  The symmetric group $S_n$ acts on
  $C^n$ by $\sigma\beta^{(j)} = \beta^{(\sigma(j))}$ and
  $\sigma\alpha_i^{(j)} = \alpha_i^{(\sigma(j))}$  and similarly on the duals for $\sigma\in S_n$, and
  this extends multiplicatively. 
\end{rem}

We have $\Sym^nC = C^n/S_n$. Since $\Sym^nC$ is smooth, we have for
$p,q\in \Z_{\ge 0}$ that
$H^q(\Sym^nC, \Omega^p_{\Sym^nC}) \cong
H^q(C^n,\Omega^p_{C^n})^{S_n}$, where the isomorphism is induced by the
quotient map $C^n \to \Sym^nC$. In the following, we identify the two
$k$-vector spaces.

\begin{rem}
	We have $H^n(C^n, \Omega_{C^n}^n) \cong k$.
	A basis is given by the element $\beta^{(1)} \cdots \beta^{(n)}$. Therefore,
        $H^n(\Sym^nC, \Omega_{\Sym^nC}^n) \cong k$ has basis
        \[\sum_{\sigma\in S_n} \beta^{(\sigma(1))} \cdots
        \beta^{(\sigma(n))} = n!\beta^{(1)}\cdots \beta^{(n)}.\]
\end{rem}

\begin{lemma}
  \label{basis-computation-lemma}
  Let $0 \le \nu\le m$ and let $I = \{i_1 < \dots < i_\nu\},J = \{j_1< \dots<
  j_\nu\} \subset \{1, \dots, g\}$ be subsets. Define
  \[
    a_{IJ} := \alpha_{i_1}^{(1)} \cdots \alpha_{i_\nu}^{(\nu)}\cdot
    \alpha_{j_1}^{\vee, (\nu+1)}\cdots \alpha_{j_\nu}^{\vee, (2\nu)}\cdot
    \beta^{(2\nu+1)}\cdots \beta^{(m+\nu)}\cdot \beta^{\vee, (m+\nu+1)} \cdots
    \beta^{\vee, (n)}
  \]
  in $H^m(C^n, \Omega_{C^n}^m)$
  and
  \[
    \alpha_{IJ} := \sum_{\sigma \in S_n} \sigma a_{IJ}.
  \]
  Then, $\alpha_{IJ}$ is an elements of $H^m(\Sym^nC,
  \Omega_{\Sym^nC}^m)$. Furthermore, $H^m(\Sym^nC, \Omega_{\Sym^nC}^m)$ has a basis given by
  $\alpha_{IJ}$, where $I,J$ runs over all same order subsets $I,J \subset \{1, \dots, g\}$ of order at most
  $m$.
\end{lemma}

\begin{proof} 
  Let $\pi_i: C^n\to C$ be the projection on the $i$'th coordinate. We 
  have
  \[\Omega_{C^n} =
    \pi_1^*\Omega_C\oplus\cdots\oplus\pi_n^*\Omega_C.\]
  This implies
  \begin{align*}
    \Omega_{C^n}^m &= \Lambda^m\Omega_{C^n}= \bigoplus_{i_1 + \cdots + i_n = m}\Lambda^{i_1}\pi_1^*\Omega_C \otimes \cdots \otimes \Lambda^{i_n}\pi_n^*\Omega_C. 
  \end{align*}
  Using the K\"unneth formula, we find
  \begin{align*} 
    H^m(C^n,\Omega_{C^n}^m) &= \bigoplus_{i_1 + \cdots + i_n = m}H^m(C^n,\Lambda^{i_1}\pi_1^*\Omega_C \otimes \cdots \otimes \Lambda^{i_n}\pi_n^*\Omega_C)\\
                            &= \bigoplus_{i_1 + \cdots + i_n = m}\bigoplus_{j_1 + \cdots + j_n = m}H^{j_1}(C,\Lambda^{i_1}\pi_1^*\Omega_C) \otimes \cdots \otimes H^{j_n}(C, \Lambda^{i_n}\pi_n^*\Omega_C).
  \end{align*} 
  We note that $H^{j}(C,\Lambda^{i}\Omega_C)$ has basis from Notation \ref{nota:basis}
  \[
    \begin{cases}
      \beta^{\vee} &\text{ if } j = i = 0,\\
      \alpha_1,\dots, \alpha_g &\text{ if } j = 0, i = 1,\\
      \alpha_1^{\vee},\dots, \alpha_g^{\vee} &\text{ if } j = 1, i = 0,\\
      \beta &\text{ if } j = i = 1.
    \end{cases}
  \]
  For $I$ and $J$ same-length tuples
  $(i_1, \dots, i_\nu)$, $(j_1, \dots j_\nu)$ with $0\le i_1 \le \dots \le i_\nu\le g$
  and $0 \le j_1 \le \dots \le j_\nu \le g$ and $\nu \le m$ define
  \[
    a_{IJ} := \alpha_{i_1}^{(1)} \cdots \alpha_{i_\nu}^{(\nu)}\cdot
    \alpha_{j_1}^{\vee, (\nu+1)}\cdots \alpha_{j_\nu}^{\vee, (2\nu)}\cdot
    \beta^{(2\nu+1)}\cdots \beta^{(m+\nu)}\cdot \beta^{\vee, (m+\nu+1)} \cdots
    \beta^{\vee, (n)}.
  \]
  (The notation $a_{IJ}$ for this is
  justified as it is a natural extension of the above definition.)
  Then, we get the basis $(\sigma a_{IJ})_{I,J,\sigma}$ for
  $H^m(C^n,\Omega_{C^n}^m)$ where $I$ and $J$ run over same-length tuples of the above form
  and $\sigma$ runs over the elements of $S_n$ not fixing $a_{IJ}$.

  Now, let $I =(i_1, \dots, i_\nu)$ and $J = (j_1, \dots, j_\nu)$ be tuples as
  above with $\nu \ge 2$. Assume that $i_{\mu} = i_{\mu+1}$ for some
  $\mu$. Then we note that $\sigma a_{IJ} = -a_{IJ}$ for $\sigma$ the
  transposition swapping $\mu$ and $\mu + 1$.

  Thus, if $v \in H^m(\Sym^nC, \Omega^m_{\Sym^nC}) = H^m(C^n,
  \Omega^m_{C^n})^{S_n}$, the $a_{IJ}$-coefficient of $v$
  vanishes for such $I,J$. Analogously, we get the same vanishing if $j_\mu = j_{\mu+1}$.
  
  In order to construct a basis of $H^m(\Sym^nC, \Omega^m_{\Sym^nC})$, we can
  therefore restrict to the $a_{IJ}$ where $i_1 < \dots < i_\nu$ and $j_1 <
  \dots < j_\nu$. Since, $S_n$ acts freely on these elements, we take sums under all $\sigma\in S_n$, to construct a
  basis of the $S_n$ fixed-points. This is the desired basis.
\end{proof} 

\begin{lemma}
  \label{multiplication-lemma}
  Let $\sigma \in S_n$ and $I,J, I', J' \subset \{1, \dots, g\}$, with $\nu = |I| =
  |J| \le m$ and $|I'| = |J'| \le m$.
  In the notation of Lemma \ref{basis-computation-lemma}, we have
  \[
    a_{IJ} \cdot \sigma a_{I'J'} = (-1)^\nu\beta^{(1)}\cdots \beta^{(n)}
  \]
  whenever $I = J', J = I'$ and $\sigma$ satisfies
  \begin{itemize}
  \item $\sigma(i) = \nu+i$ and $\sigma(\nu+i) = i$ for $1 \le i \le \nu$ and
  \item $\sigma(\{2\nu+1, \dots, m+ \nu\}) = \{m+ \nu+ 1, \dots, n\}$.
  \end{itemize}
  Otherwise
  $a_{IJ} \cdot \sigma a_{I'J'} = 0$.

  \begin{proof}
    The last assertion follows from the multiplication table in Remark
    \ref{rem-multiplication-on-curve}.

    Assume that $I = J', J = I'$ and that $\sigma$ satisfies the condition of
    the lemma. Then, we have
    \[
      \sigma a_{JI} = \alpha_{j_1}^{(\nu + 1)} \cdots \alpha_{j_\nu}^{(2\nu)}\cdot
      \alpha_{i_1}^{\vee, (1)}\cdots \alpha_{i_\nu}^{\vee, (\nu)}\cdot
      \beta^{(\sigma(2\nu+1))}\cdots \beta^{(\sigma(m+\nu))}\cdot \beta^{\vee, (\sigma(m+\nu+1))} \cdots
      \beta^{\vee, (\sigma(n))}.
    \]
    Since the degrees of $\beta^{(i)}$ and $\beta^{\vee,(i)}$ are even, we can
    rearrange these without introducing a sign. Thus, the above equation
    can be rewritten as
    \[
      \sigma a_{JI} = \alpha_{j_1}^{(\nu + 1)} \cdots \alpha_{j_\nu}^{(2\nu)}\cdot
      \alpha_{i_1}^{\vee, (1)}\cdots \alpha_{i_\nu}^{\vee, (\nu)}\cdot
      \beta^{\vee, (2\nu+1)}\cdots \beta^{\vee, (m+\nu)}\cdot \beta^{(m+\nu+1)} \cdots
          \beta^{(n)}.
    \]
    Define $B' := \beta^{\vee, (2\nu+1)}\cdots \beta^{\vee, (m+\nu)}\cdot \beta^{(m+\nu+1)} \cdots
    \beta^{(n)}$ and,
    similarly, define $B := \beta^{(2\nu+1)}\cdots \beta^{(m+\nu)}\cdot
    \beta^{\vee, (m+\nu+1)} \cdots \beta^{\vee, (n)}$.
    Then, since $B$ and $B'$ have even degree, we have that
    \begin{align*}
      a_{IJ} \cdot \sigma a_{JI}
      &= \alpha_{i_1}^{(1)} \cdots \alpha_{i_\nu}^{(\nu)}\cdot
    \alpha_{j_1}^{\vee, (\nu+1)}\cdots \alpha_{j_\nu}^{\vee, (2\nu)} \cdot B
        \cdot \alpha_{j_1}^{(\nu + 1)} \cdots \alpha_{j_\nu}^{(2\nu)}\cdot
        \alpha_{i_1}^{\vee, (1)}\cdots \alpha_{i_\nu}^{\vee, (\nu)}\cdot B'\\
      &= \alpha_{i_1}^{(1)} \cdots \alpha_{i_\nu}^{(\nu)}\cdot
    \alpha_{j_1}^{\vee, (\nu+1)}\cdots \alpha_{j_\nu}^{\vee, (2\nu)}
        \cdot \alpha_{j_1}^{(\nu + 1)} \cdots \alpha_{j_\nu}^{(2\nu)}\cdot
        \alpha_{i_1}^{\vee, (1)}\cdots \alpha_{i_\nu}^{\vee, (\nu)}\cdot B\cdot B'
    \end{align*}
    All $\alpha^{(i)}$'s and $\alpha^{\vee, (i)}$'s have odd degree. Thus, if we rearrange these elements by a permutation
    $\tau$, we have to multiply with $(-1)^{\sgn(\tau)}$. Therefore, if we
    rearrange
    \[
      \alpha_{i_1}^{(1)} \cdots \alpha_{i_\nu}^{(\nu)}\cdot
      \alpha_{j_1}^{\vee, (\nu+1)}\cdots \alpha_{j_\nu}^{\vee, (2\nu)}
      \cdot \alpha_{j_1}^{(\nu + 1)} \cdots \alpha_{j_\nu}^{(2\nu)}\cdot
      \alpha_{i_1}^{\vee, (1)}\cdots \alpha_{i_\nu}^{\vee, (\nu)}\cdot
      B\cdot B'
    \]
    as
    \[
      \alpha_{i_1}^{(1)} \cdots \alpha_{i_\nu}^{(\nu)}\cdot
      \alpha_{j_1}^{\vee, (\nu+1)}\cdots \alpha_{j_\nu}^{\vee, (2\nu)}
      \cdot
      \alpha_{i_1}^{\vee, (1)}\cdots \alpha_{i_\nu}^{\vee, (\nu)}\cdot
      \alpha_{j_1}^{(\nu + 1)} \cdots \alpha_{j_\nu}^{(2\nu)}\cdot
      B\cdot B',
    \]
    we are swapping $\nu$ times, so this introduces the sign $(-1)^\nu$. If we
    next rearrange this into 
    \[
      \alpha_{i_1}^{(1)}\alpha^{\vee,(1)}_{i_1} \cdots
      \alpha_{i_\nu}^{(\nu)}\alpha_{i_\nu}^{\vee, (\nu)}\cdot
      \alpha_{j_1}^{\vee, (\nu+1)}\alpha^{(\nu+1)}_{j_1}\cdots
      \alpha_{j_\nu}^{\vee, (2\nu)}\alpha_{j_\nu}^{(2\nu)}\cdot
      B\cdot B',
    \]
    we are moving $\alpha_{i_1}^{\vee, (1)}$ by $2\nu-1$ elements to the
    left, then we are moving $\alpha_{i_2}^{\vee, (2)}$ by $2\nu-2$ elements to the
    left, $\alpha_{i_3}^{\vee, (3)}$ by $2\nu-3$ elements, and so on. Therefore,
    the resulting permutation has sign $\sum_{i=1}^{2\nu} 2\nu -i =
    \sum_{i=0}^{2\nu-1} i = \frac{2\nu(2\nu-1)}{2} = \nu(2\nu -1)$. Thus, we
    have
    \begin{align*}
      a_{IJ} \cdot \sigma a_{JI}
      &= \underset{=1}{\underbrace{(-1)^{\nu}(-1)^{\nu(2\nu-1)}}}\cdot
        \alpha_{i_1}^{(1)}\alpha^{\vee,(1)}_{i_1} \cdots
        \alpha_{i_\nu}^{(\nu)}\alpha_{i_\nu}^{\vee, (\nu)}\cdot
        \alpha_{j_1}^{\vee, (\nu+1)}\alpha^{(\nu+1)}_{j_1}\cdots
        \alpha_{j_\nu}^{\vee, (2\nu)}\alpha_{j_\nu}^{(2\nu)}\cdot
        B\cdot B'\\
      &= \underset{=\beta^{(1)}}{\underbrace{\alpha_{i_1}^{(1)}\alpha^{\vee,(1)}_{i_1}}} \cdots
        \underset{=\beta^{(\nu)}}{\underbrace{\alpha_{i_\nu}^{(\nu)}\alpha_{i_\nu}^{\vee, (\nu)}}}\cdot
        \underset{=-\beta^{(\nu+1)}}{\underbrace{\alpha_{j_1}^{\vee, (\nu+1)}\alpha^{(\nu+1)}_{j_1}}}\cdots
        \underset{=-\beta^{(2\nu)}}{\underbrace{\alpha_{j_\nu}^{\vee, (2\nu)}\alpha_{j_\nu}^{(2\nu)}}}\cdot
        B\cdot B'\\
      &= (-1)^\nu \beta^{(1)}\cdots \beta^{(2\nu)} \underset{=
        \beta^{(2\nu+1)}\cdots \beta^{(n)}}{\underbrace{B\cdot B'}}\\
      &= (-1)^\nu \beta^{(1)}\cdots \beta^{(n)}.\qedhere
    \end{align*}
  \end{proof}
\end{lemma}

\begin{lemma}
  \label{lem:trace-finite-change}
  Let $f\colon X \to Y$ be a finite morphism of degree $d$ between smooth,
  projective, $n$-dimensional schemes $X$ and $Y$
  over $k$. Then, we have the commutative diagram
  \[
    \begin{tikzcd}
      H^n(X, \Omega^n_X) \ar[r, "\Tr_{X/k}"] & k.\\
      H^n(Y, \Omega^n_X) \ar[ur, "d\cdot \Tr_{Y/k}"']\ar[u,
      "f^\ast"] &
    \end{tikzcd}
  \]

  \begin{proof}
    This argument was explained to the first named author by Marc Levine. Let
    $f_*: H^n(X,\Omega^n_X)\to H^n(Y,\Omega^n_Y)$ be the pushforward as
    constructed by Srinivas \cite{Srinivas}. Then by \cite[Theorem
    1]{Srinivas}, the diagram
    \[
    \begin{tikzcd}
      H^n(X, \Omega^n_X) \ar[r, "\Tr_{X/k}"]\ar[d, "f_\ast"] & k\\
      H^n(Y, \Omega^n_Y) \ar[ur, "\Tr_{Y/k}"'] &
    \end{tikzcd}
  \]
  commutes and we have $f_\ast f^\ast = d\cdot \operatorname{id}$. This yields the desired result. 
  \end{proof}
\end{lemma}

\begin{proof}[Proof of Theorem \ref{main theorem} in case $n=2m$ is even and
  $k$ has characteristic zero]
  We will apply the Motivic Gauss-Bonnet Theorem, and thus, we need to compute the form 
  \begin{equation}\label{form q}
  \begin{tikzcd} [row sep= small]
    H^m(\Sym^nC, \Omega_{\Sym^nC}^m)\times H^m(\Sym^nC,
    \Omega_{\Sym^nC}^m) \arrow[d, "\cup"]\\
    H^n(\Sym^nC,
    \Omega_{\Sym^nC}^n) \arrow[d, "\Tr"]\\
    k.
   \end{tikzcd} 
  \end{equation}
  By Lemma \ref{basis-computation-lemma}, the vector space $H^m(\Sym^nC, \Omega_{\Sym^nC}^m)$ has basis $\alpha_{IJ}$ for
  $I,J$ running over the same size subsets $I,J
  \subset \{1, \dots, g\}$ with $|I| = |J| \le m$. By Lemma
  \ref{multiplication-lemma}, we know that $\alpha_{IJ}\alpha_{I'J'} = 0$
  whenever $I \ne J'$ or $J \ne I'$. Thus, it remains to compute
  the cup product $\alpha_{IJ}\alpha_{JI}$. For this, note that
  \[
    \alpha_{IJ}\alpha_{JI} = \sum_{\sigma,\tau \in S_n} (\sigma a_{IJ})\cdot
    (\tau a_{JI})= \sum_{\sigma, \tau \in S_n} \sigma (a_{IJ}\cdot
    (\sigma^{-1}\tau)a_{JI}) = \sum_{\sigma, \tau \in S_n} \sigma
    (a_{IJ}\cdot \tau a_{JI})
  \]
  where in the last step, we replace $\sigma^{-1}\tau$ by $\tau$. Since $a_{IJ}\cdot \tau a_{JI}$ is a multiple of $\beta^{(1)}\cdots
  \beta^{(n)}$ and $S_n$ acts trivially on this element, we get by Lemma \ref{multiplication-lemma}
  \begin{align*} 
    \alpha_{IJ}\alpha_{JI} &= \sum_{\sigma, \tau \in S_n} \sigma
    (a_{IJ}\cdot \tau a_{JI})\\
    &= n!\cdot \sum_{\tau \in S_n}
    a_{IJ}\cdot \tau a_{JI} \\
    &= n!\cdot \sum_{\tau \in S_n'}
    (-1)^{|I|}\beta^{(1)}\cdots \beta^{(n)} \\
    &= (-1)^{|I|}\cdot n!\cdot
    |S_n'|\cdot \beta^{(1)}\cdots \beta^{(n)}
  \end{align*}
  where $S_n' \subset S_n$ is the subset of all permutations satisfying the
  non-vanishing condition of Lemma \ref{multiplication-lemma} for $I$ and
  $J$. Thus, in order to compute the product, we have to determine the cardinality
  of $S_n'$. Note that
  \[
    S_n' = \{\tau\sigma\mid \sigma \in S_{m-|I|}\times
    S_{m-|I|}\},
  \]
  where we embed $S_{m-|I|}\times
  S_{m-|I}$ in $S_n$ by letting it permute the last $2(m-|I|)$ elements of
  $\{1, \dots, n\}$ and
  where $\tau\in S_n$ is the permutation defined as follows: We set $\tau(i) = |I|+i$ and $\tau(|I|+i) = i$ for $1 \le i \le |I|$,
  and $\tau(2|I|+j) = m+ |I|+ j$ and $\tau(m+|I|+j) = 2|I| + j$ for $1 \le
  j \le m-|I|$.
  Thus, $|S_n'| = |S_{m-|I|}\times S_{m-|I|}| =
  ((m-|I|)!)^2$, and we have
  \[
    \alpha_{IJ}\alpha_{JI} = (-1)^{|I|}\cdot n!\cdot
    ((m-|I|)!)^2\cdot \beta^{(1)}\cdots \beta^{(n)}.
  \]
  Since Lemma
  \ref{lem:trace-finite-change} yields $\Tr(\beta^{(1)}\cdots \beta^{(n)}) = (n!)^{-1}$, we obtain
  \[
    \Tr(\alpha_{IJ}\alpha_{I'J'}) =
    \begin{cases}
      (-1)^{|I|}\cdot ( (m-|I|)!)^2 & \text{if }I=J' \text{ and } J =
                                             I'\\
      0 & \text{otherwise.}
    \end{cases}
  \]
  Thus, the matrix representing the bilinear form with respect to the basis
  $(\alpha_{IJ})$ takes the form of a diagonal block matrix with blocks
  $A_{\{I,J\}}$ for $I,J \subset \{1, \dots, g\}$ same-size subsets of order at
  most $m$. Each block $A_{\{I,J\}}$ takes the form
  \[
    A_{\{I,J\}} =
    \begin{pmatrix}
      (-1)^{|I|}\cdot ((m-|I|)!)^2
    \end{pmatrix}
  \]
  if $I = J$ and
  \[
    A_{\{I,J\}} =
    \begin{pmatrix}
      0 & (-1)^{|I|}\cdot ( (m-|I|)!)^2\\
      (-1)^{|I|}\cdot ( (m-|I|)!)^2 & 0
    \end{pmatrix}
  \]
  if $I\neq J$. The matrix $A_{\{I,I\}}$ represents the form $\langle (-1)^{|I|}\cdot
  ( (m-|I|)!)^2\rangle = \langle (-1)^{|I|}\rangle \in \GW(k)$, and the
  matrix $A_{\{I,J\}}$ with $I \ne J$ represents the hyperbolic form in
  $\GW(k)$. Thus, we have that the quadratic form $q$ in \eqref{form q} is given by
  \[
    q=\sum_{I \subset \{1, \dots, g\}} \langle (-1)^{|I|}\rangle + l\cdot H
  \]
  for some $l \in \Z$ where the sum only runs over subsets of order at most
  $m$. Since $\{1, \dots, g\}$ contains $\binom g i$ subsets of order $i$, we
  get
  \[
    q = \sum_{i=0}^m\binom gi \langle (-1)^i\rangle + lH.
  \]
  Combining this computation with Remark \ref{rem:rank} yields the first
  formula of the theorem. For the second formula, we plug in $H = \langle 1\rangle + \langle -1\rangle$ and then note that
  \[
    \sum_{i=0}^m\binom{g}{i}(-1)^i = (-1)^m\binom{g-1}{m}.\qedhere
  \]
\end{proof}

\section{Compatibility with the power structure}
Recall that a power structure on a commutative ring $R$ is a map
$(1+tR[[t]])\times R\to 1+tR[[t]], (f(t),r)\mapsto f(t)^r$ satisfying certain axioms, see for example
\cite[Section 2]{Pajwani-PalPS} for a precise definition. Under some
finiteness assumptions, it suffices to define $(1-t)^{-r} = \sum_{n=0}^\infty
a_n(r)t^n$ for $r \in R$, see \cite[Proposition 1]{Gusein-Zade-Luengo-Melle-HernandezHS}. That is, a power structure can be given by functions $a_n: R\to R$ for $n\in\mathbb{Z}_{\geq 0}$ such that:
\begin{itemize}
	\item $a_0 = 1$,
	\item $a_1$ is the identity,
	\item $a_n(0) = 0$ for all $n\geq 1$,
	\item $a_n(1) = 1$ for all $n\geq 1$, and
	\item $a_n(r+s) = \sum_{i=0}^na_i(r)a_{n-i}(s)$ for all $r,s\in R$.
\end{itemize} 
A ring homomorphism $\phi: R\to R'$ between rings that both carry a power structure is said to respect the power structures if $\phi(a_n^R(r)) = a_n^{R'}(\phi(r))$ for all $r\in R$ and $n\in\mathbb{Z}_{\geq 0}$. 

For the remainder of this section, we assume that $k$ is of characteristic zero.

There is a power structure on $K_0(\Var_k)$ given by $a_n([X]) = [\Sym^nX]$. This was proven by Gusein-Zade, Luengo and Melle-Hern\'andez  \cite{Gusein-Zade-Luengo-Melle-HernandezPS}. 

Pajwani and P\'al \cite{Pajwani-PalPS} constructed a power structure on
$\GW(k)$ given by
\[
  a_n(\langle a \rangle) = \langle a^n\rangle + \frac{n(n-1)}{2}(\langle a
  \rangle + \langle 2 \rangle - \langle 1 \rangle - \langle 2a \rangle).
\]
\begin{rem} 
  This power structure agrees with McGarraghy's non-factorial power structure
  \cite{McGarraghySP} on $\GW(k)$ on the subgroup generated by all $\langle a
  \rangle$ satisfying $a_n(\langle a \rangle) = \langle a^n\rangle$, see
  \cite[Lemma 2.9, Lemma 3.27 and Corollary 3.28]{Pajwani-PalPS}. The quadratic
  forms $\langle 1 \rangle$ and $\langle -1 \rangle$ are in this subgroup. 
\end{rem} 
Pajwani and P\'al prove that $\chi_c: K_0(\Var_k)\to \GW(k)$ respects the power structures for
zero dimensional schemes, i.e. $\chi_c(\Sym^nX/k) = a_n(\chi_c(X/k))$ for $X$ a
zero dimensional scheme. This was extended to a larger class of schemes by
Pajwani, Rohrbach, and the second named author in \cite{PajwaniRohrbachViergever}, who also showed
compatibility for curves of genus $0$ and $1$. We now extend this result to
all curves, i.e. we check that all curves are symmetrizable in the sense of
\cite[Definition 4.1]{PajwaniRohrbachViergever}. 

For $C$ a smooth, projective curve of genus $g$, we have $\chi(C/k) =
(1-g)H$. We therefore need to evaluate $a_n$ on negative multiples of
$H$. Note that McGarraghy \cite[Corollary 4.13 and 4.14]{McGarraghySP} gave a
complete description of what happens for positive multiples.

\begin{lemma}
	Let $l\in\mathbb{Z}_{\geq 1}$. Then $a_n(-lH) = (-1)^n
        \sum_{i=0}^n\binom{l}{i}\binom{l}{n-i}\langle -1
        \rangle^{n-i}$.
\end{lemma} 

\begin{proof}
	One can prove by induction on $l$ that 
	\[
		a_n(-l\langle 1 \rangle ) = (-1)^n\binom{l}{n}\langle 1 \rangle
	\text{ and }
		a_n(-l\langle -1 \rangle ) = (-1)^n\binom{l}{n}\langle -1 \rangle^n.\] 
	We obtain 
	\begin{align*}
		a_n(-lH) &= \sum_{i=0}^na_i(-l\langle 1 \rangle)a_{n-i}(-l\langle -1 \rangle) \\
		&= (-1)^n\sum_{i=0}^n\binom{l}{i}\binom{l}{n-i}\langle -1\rangle^{n-i} 
	\end{align*}
	as desired. 
\end{proof}

\begin{rem}\label{remark: same rank}
	We note that for $l,n\in\mathbb{Z}_{\geq 1}$, the rank of $a_n(-lH)$ equals 
	$$(-1)^n\sum_{i=0}^n\binom{l}{i}\binom{l}{n-i} = (-1)^n\binom{2l}{n}$$
	by the Chu-Vandermonde identity. Let $C$ be a smooth projective curve over $k$ of genus $g$, then we know from Theorem \ref{main theorem} that $\text{rank}(\chi(\Sym^nC/k)) = (-1)^n\binom{2g-2}{n}$. Plugging in $l = g-1$ to the above formula, we find that $\text{rank}(a_n(1-g)H) = (-1)^n\binom{2g-2}{n}$. Therefore, we see that the ranks of $\chi(\Sym^nC/k)$ and $a_n(\chi(C/k))$ agree.
\end{rem}

There is another way to compute $a_n(-lH)$, suggested by Marc Levine. 

\begin{lemma}\label{lemma: am as coefficient of polynomial}
	Let $l\in\mathbb{Z}_{\geq 1}$. Then $a_n(-lH)$ is the coefficient of $t^n$ in the polynomial $ (\langle 1 \rangle -Ht+\langle -1 \rangle t^2)^l $.
\end{lemma} 

\begin{proof}
	For $\langle a \rangle\in\GW(k)$, let $a_t(\langle a \rangle) :=
        (\langle 1 \rangle-\langle 1\rangle t)^{-\langle a \rangle} = \sum_{i=0}^\infty a_i(\langle a
        \rangle)t^i$. Then $a_t(\langle a \rangle + \langle b \rangle) =
        a_t(\langle a\rangle)a_t(\langle b\rangle)$. Now, we have $a_t(-\langle 1 \rangle) = \langle 1 \rangle-\langle 1 \rangle t$
	and also 
	$a_t(-\langle -1 \rangle) = \langle 1 \rangle -\langle -1 \rangle t.$
	This yields
	\begin{align*}
		a_t(-lH) &= ((\langle 1 \rangle -\langle 1 \rangle t)(\langle 1 \rangle -\langle -1\rangle t))^l \\
		&= (\langle 1 \rangle -Ht+\langle -1 \rangle t^2)^l 
	\end{align*} 
	as desired. 
\end{proof}

\begin{prop}\label{proposition: ps respected smooth projective curves}
	Let $C$ be a smooth, projective curve of genus $g$ and let
        $n\in\mathbb{Z}_{\geq 0}$. Then $$\chi(\Sym^nC/k)=a_n(\chi(C/k)).$$
\end{prop} 

\begin{proof}
	The two forms have the same rank by Remark \ref{remark: same
          rank}. Given that they are both hyperbolic for $n$ odd, this yields
        the result in that case. Thus, assume that $n = 2m$ is even from now on. 
	 
	Note that both forms only consist of terms $\langle 1 \rangle$ and $\langle -1 \rangle$. Therefore, given that we know the ranks are the same, the two are the same if for both forms the difference between the number of $\langle 1\rangle$'s and $\langle -1\rangle$'s is the same. Note that if we work over the real numbers this difference is the signature of the quadratic form. By Theorem \ref{main theorem}, this difference is 
	$(-1)^m\binom{g-1}{m}$ for $\chi(\Sym^{2m}C/k)$.

	On the other hand, denote this difference for $a_n((1-g)H)$ by $s$. Lemma \ref{lemma: am as coefficient of polynomial} yields that $a_n((1-g)H)$ is the coefficient of $t^n$ in 
	\begin{align*}
			(\langle 1 \rangle -Ht+\langle -1 \rangle t^2)^{g-1} &= \sum_{i=0}^{g-1} \binom{g-1}{i} (\langle 1 \rangle + \langle -1 \rangle t^2)^i (-Ht)^{g-1-i}.
	\end{align*}
	All summands with $i\neq g-1$ are a multiple of $H$ and thus do not contribute to~$s$. Thus, $s$ is completely determined by the degree $2m$-term of 
	\begin{align*} 
		(\langle 1\rangle + \langle -1\rangle t^2)^{g-1} &= \sum_{i=0}^{g-1} \binom{g-1}{i}\langle -1 \rangle^{i}t^{2i},
	\end{align*} 
	which is $\binom{g-1}{m}\langle -1\rangle^mt^{2m}$. This shows $s = (-1)^m\binom{g-1}{m}$, which proves the proposition. 
\end{proof}

We can extend this to show that the power structure is respected for all curves. 

\begin{theorem}\label{theorem power structures respected}
	Let $C$ be a curve and let
        $n\in\mathbb{Z}_{\geq 0}$. Then $$\chi_c(\Sym^nC/k) =
        a_n(\chi_c(C/k)).$$
\end{theorem}

\begin{proof}
	Let $\tilde{C}$ be the normalization of $C$. Let $S\subset C$ be the
        locus of singular points on $C$ and let $\tilde{S}\subset\tilde{C}$ be
        the points lying over $S$. Then, the normalization map $\pi: \tilde{C}\to C$ is an isomorphism outside of the singular points of $C$ and so 
	$$[C]  = [\tilde{C}] - [\tilde{S}] + [S]$$
	in $\tilde{K_0}(\text{Var}_k)$. By Proposition \ref{proposition: ps respected smooth projective curves}, the power structure is respected on $\tilde{C}$. Also, as $S$ and $\tilde{S}$ are zero dimensional, the power structure is respected on those by \cite[Theorem 4.1]{Pajwani-PalPS}. By \cite[Lemma 2.9]{Pajwani-PalPS}, the set of elements of $\tilde{K_0}(\text{Var}_k)$ for which the power structure is respected forms a subgroup, implying that the power structure is respected for $C$.
\end{proof}

\phantomsection 
\addcontentsline{toc}{section}{Bibliography}
{\small
\bibliographystyle{plain}
\bibliography{Sources}}

~\\
	
%	\newpage 
    \noindent Lukas F. Br\"oring \\
    Universit\"at Duisburg-Essen \\
    Fakult\"at f\"ur Mathematik, Thea-Leymann-Str. 9, 45127 Essen, Germany \\
	E-Mail: \href{mailto:lukas.broering@uni-due.de}{ lukas.broering@uni-due.de}\\

	\noindent Anna M. Viergever \\
	Leibniz Universit\"at Hannover\\
	Fakult\"at f\"ur Mathematik, Welfengarten 1, 30167 Hannover, Germany \\
	E-Mail: \href{mailto:viergever@math.uni-hannover.de}{ viergever@math.uni-hannover.de} \\ \\ 

	\noindent Keywords: motivic homotopy theory, refined enumerative geometry, other fields\\ 
	Mathematics Subject Classification: 14G27, 14N10, 14F42

\end{document}